\documentclass[preprint,12pt,authoryear]{elsarticle}
\usepackage{amsmath}
\usepackage{amssymb}
\usepackage{amsthm}

\newtheorem{theorem}{\indent Theorem}
\newtheorem{lemma}{\indent Lemma}

\newtheorem{corollary}{Corollary}
\newtheorem{definition}{\indent Definition}
\newtheorem{remark}{\indent Remark}

\providecommand{\abs}[1]{\lvert#1\rvert}
\providecommand{\norm}[1]{\displaystyle\left\lVert#1\right\rVert}

\journal{Journal of Computational and Applied Mathematics} 

\begin{document}

\begin{frontmatter}



\title{Asymptotical properties of social network dynamics on time scales}
\author{Aleksey Ogulenko}
\ead{ogulenko.a.p@onu.edu.ua}
\address{I. I. Mechnikov Odesa National University, Dvoryanska str, 2, Odesa, Ukraine, 65082}


\begin{abstract}
	In this paper we develop conditions for various types of stability in
	social networks governed by {\it Imitation of Success} principle.
	Considering so-called Prisoner's Dilemma as the base of node-to-node
	game in the network we obtain well-known Hopfield neural network model.
	Asymptotic behavior of the original  model and dynamic Hopfield model
	have a certain correspondence. To obtain more general results, we
	consider Hopfield model dynamic system on time scales.  Developed
	stability conditions combine main parameters of network structure such
	as network size and maximum relative nodes' degree with the main
	characteristics of time scale, nodes' inertia and resistance, rate of
	input-output response. 
\end{abstract}


\end{frontmatter}

\section{Introduction}
A social network is the set of people or groups of people with some pattern of
links or interconnection between them. Processes taking place on social
networks often may be interpreted as information transition.

The aim of this paper is to consider asymptotic properties of collective
opinion formation in social networks with general topology. Transition of
opinion between linked nodes will be modelled by game-theoretical mechanism.
Total payoff may be a key factor to choose one of the two alternative
strategies, cooperation or defection in opinion propagation.  Such type of
dynamics is called {\it Imitation of Success}.  An opposite (in some sense)
kind of model is for example the {\it Voter Model}.  Last named model assume
imitation of a behavior of uniformly random chosen neighbor node and games
payoff has no affects on state updating of the particular node.
In paper \cite{monshadi2011} was considered model called {\it Weak 
Imitation of Success}. This updating rule is mixture of IS and VM rules:
dependently of some parameter $\varepsilon$ behavior of updating node
may be close to one of the two types of dynamics.

The analysis of the total payoff function for so-called Prisoner's Dilemma
leads us to well-known Hopfield neural network model \cite{hopfield1984}.
Asymptotic behavior of the direct node-to-node model and dynamic Hopfield
model have a certain correspondence. To obtain more general results, we
consider Hopfield model on time scale. This problem is discussed in detail in
\cite{martynyuk2012-green}, but we develop more direct and precise conditions
for stability of
the social network behavior.

\section{Preliminary results.}

We now present some basic information about time scales according to \cite{bohner2012dynamic}.
A time scale is defined as a nonempty closed subset of the set of
real numbers and denoted by $\mathbb{T}$.  The properties of the time scale are
determined by the following three functions:
\begin{itemize}
	\item[(i)] the forward-jump operator: $\sigma(t) = \inf \left\{s \in \mathbb{T}: s > t\right\}$;
	\item[(ii)] the backward-jump operator: $\rho(t) = \sup \left\{s \in \mathbb{T}: s < t \right\}$ 
		(in this case, we set $\inf\varnothing  = \sup\mathbb{T}$ and $\sup\varnothing = \inf\mathbb{T}$);
	\item[(iii)] the granularity function $\mu(t) = \sigma(t) - t$.
\end{itemize}

The behavior of the forward- and backward-jump operators at a given point of
the time scale specifies the type of this point. The corresponding
classification of points of the time scale is presented in Table 1.
\begin{table}[h]
\label{tab:pointclass}
    \caption{Classification of time scale's points}
    \begin{center}
        \begin{tabular}{|l|l|}
            \hline
            $t$ right-scattered & $t < \sigma(t)$ \\
            $t$ right-dense    & $t = \sigma(t)$ \\
            $t$ left-scattered  & $ \rho(t) < t $ \\
            $t$ left-dense     & $ \rho(t) = t $ \\
            $t$ isolated      & $\rho(t) < t < \sigma(t)$ \\
            $t$ dense           & $\rho(t) = t = \sigma(t)$ \\
            \hline
        \end{tabular}
    \end{center}
\end{table}

We define a set $\mathbb{T}^\kappa$ in the following way:
\[
    \mathbb{T}^\kappa =
        \begin{cases}
            \mathbb{T}\setminus \left\{M\right\}, & \text{ if }
            \exists \text{ right scattered point } M \in \mathbb{T}: M = \sup\mathbb{T}, \sup\mathbb{T}<\infty  \\
            \mathbb{T} , & \text{ otherwise}.
        \end{cases}
\]
In what follows, we set $\left[a, b\right] = \left\{t \in \mathbb{T}: a \leqslant t \leqslant b \right\}$.

\begin{definition}
    Let $f:\mathbb{T} \rightarrow \mathbb{R}$ and $t \in \mathbb{T}^\kappa$.
    The number $f^{\Delta}(t)$ is called $\Delta$-derivative of function $f$ at the point $t$,
    if $\forall \varepsilon > 0$ there exists a neighborhood $U$ of the point $t$ 
    (i.~e., $U = (t - \delta, t + \delta) \cap \mathbb{T}, \delta < 0$) such that 
    \[
    \abs{f(\sigma(t)) - f(s) - f^{\Delta}(t)(\sigma(t)-s)} \leqslant
        \varepsilon\abs{\sigma(t) - s} \quad \forall s \in U.
    \]
\end{definition}

\begin{definition}
    If $f^{\Delta}(t)$ exists $\forall t \in \mathbb{T}^\kappa$,
    then $f:\mathbb{T} \rightarrow \mathbb{R}$ is called $\Delta$-differentiable on 
    $\mathbb{T}^\kappa$. The function $f^{\Delta}(t): \mathbb{T}^\kappa \rightarrow \mathbb{R}$
    is called the delta-derivative of a function $f$ on  $\mathbb{T}^\kappa$. 
\end{definition}
If $f$ is differentiable with respect to $t$ then $f(\sigma(t)) = f(t) + \mu(t)f^{\Delta}(t)$.
\begin{definition}
	The function $f:\mathbb{T} \rightarrow \mathbb{R}$ is called regular if
	it has finite right limits at all right-dense points of the time scale $\mathbb{T}$ 
	and finite left limits at all points left-dense points of~$\mathbb{T}$.
\end{definition}

\begin{definition}
	The function $f:\mathbb{T} \rightarrow \mathbb{R}$ is called
	$rd$-continuous if it is continuous at the right-dense points and has
	finite left limits at the left-dense points.  The set of these
	functions is denoted by $C_{rd} = C_{rd}(\mathbb{T}) = 
	C_{rd}(\mathbb{T}; \mathbb{R})$.
\end{definition}

The indefinite integral on the time scale takes the form
\[
    \int{f(t) \Delta t} = F(t) + C,
\]
where $C$ is integration constant and $F(t)$ is the preprimitive for $f(t)$.
If the relation $F^{\Delta}(t) = f(t)$ where $f:\mathbb{T}\rightarrow\mathbb{R}$ 
is an $rd$-continuous function, is true for all $t \in \mathbb{T}^\kappa$ then 
$F(t)$ is called the primitive of the function $f(t)$. If $t_0 \in \mathbb{T}$ then
$F(t) = \int\limits_{t_0}^t{f(s) \Delta s}$ for all $t$.
For all $r,s \in \mathbb{T}$ the definite $\Delta$-integral is defined as follows:
\[
    \int\limits_r^s{f(t)\Delta t} = F(s) - F(r).
\]

\begin{definition}
	For any regular function $f(t)$ there exists a function $F$
	differentiable in the domain $D$ and such that the equality
	$F^{\Delta}(t) = f(t)$ holds for all $t \in D$. This function is
	defined ambiguously. It is called the preprimitive of $f(t)$.
\end{definition}

\begin{definition}
A function $p:\mathbb{T} \rightarrow \mathbb{R}$ is called regressive (positive regressive) if
    \[
        1 + \mu(t)p(t) \neq 0, \quad (1 + \mu(t)p(t) > 0), \qquad t \in \mathbb{T}^\kappa.
    \]
The set of regressive (positive regressive) and $rd$-continuous functions is denoted by 
$\mathcal{R} = \mathcal{R}(\mathbb{T})$ ($\mathcal{R}^+ = \mathcal{R}^+(\mathbb{T})$).
\end{definition}

\begin{definition}
For any $p, q \in \mathcal{R}$ by definition put  
    \[
        (p \oplus q)(t) = p(t) + q(t) + p(t)q(t)\mu(t), \quad t \in \mathbb{T}^\kappa.
    \]
\end{definition}
It is easy to see that the pair $\left(\mathcal{R}, \oplus\right)$ is an
Abelian group.  As shown in [Bohner Peterson], a function $p$ from the class
$\mathcal{R}$ can be associated with a function $e_p(t, t_0)$ which is the
unique solution of Cauchy problem
\[
    y^{\Delta} = p(t)y, \quad y(t_0) = 1.
\]
The function $e_p(t,t_0)$ is an analog, by its properties, of the exponential
function defined on~$\mathbb{R}$.

Let us consider dynamic system on time scale $\mathbb{T}$:
\begin{equation}
	\label{general-equation}
	\begin{cases}
		x^\Delta = f(t, x)\\
		x(t_0) = x_0.
	\end{cases}
\end{equation}
In following formulations we denote solution of \eqref{general-equation} by $x(t; t_0, x_0)$.
\begin{definition}
	The equilibrium state $x = x^*$ of the system \eqref{general-equation}
	is uniformly stable if $\forall \varepsilon > 0$ there exists $\delta =
	\delta(\varepsilon)$, such that 
	\[
		\norm{x_0 - x^*} < \delta \implies 
			\norm{x(t; t_0, x_0) - x^*} < \varepsilon, 
		\qquad \forall t \in [t_0, +\infty]_{\mathbb{T}}, t_0 \in \mathbb{T}.
	\]
\end{definition}

\begin{definition}
	The equilibrium state $x = x^*$ of the system \eqref{general-equation}
	is uniformly asymptotically stable if it is uniformly stable 
	and there exists $\Delta > 0$ such that 
	\[
		\norm{x_0 - x^*} < \Delta \implies 
		\lim\limits_{t \to +\infty}{\norm{x(t; t_0, x_0) - x^*}} = 0, 
		\qquad \forall t_0 \in \mathbb{T}.
	\]
\end{definition}

\begin{definition}
	The equilibrium state $x = x^*$ of the system \eqref{general-equation}
	is uniformly exponentially stable if there exist constants $\alpha, \beta > 0$ 
	($\beta \in \mathcal{R}^+$) such that 
	\[
		\norm{x(t, t_0, x_0)} \leqslant \norm{x(t_0)}\alpha e_{-\beta}(t, t_0),  
			\quad t \geqslant t_0 \in \mathbb{T},
	\]
	for all $t_0 \in \mathbb{T}$ and $x(t_0) \in \mathbb{R}^n$.
\end{definition}

In further definitions and theorems \ref{lyapunov-theorem-1} --
\ref{lyapunov-theorem-3}  we assume $f(t, 0) = 0$ for all $t \in \mathbb{T}, t
\geqslant t_0$ and $x_0 = 0$ so that $x = 0$ is a solution to
equation~\eqref{general-equation}. For more details see \cite{hoffacker2005}.

\begin{definition}
	A function $\psi : \left[0, r\right] \to \left[0, \infty\right)$ 
	is of class $\mathcal{K}$ if it is well-defined, continuous, and
	strictly increasing on $\left[0, r\right]$ with $\psi(0) = 0$.
\end{definition}

\begin{definition}
	A continuous function $P : \mathbb{R}^n \to \mathbb{R}$ with $P(0) = 0$ 
	is called positive definite (negative definite) on $D$ if there exists a 
	function $\psi \in \mathcal{K}$, such that $\psi\left(\|x\|\right) \leqslant P(x)$
	($\psi\left(\|x\|\right) \leqslant -P(x)$) for $x \in D$.
\end{definition}

\begin{definition}
	A continuous function $P : \mathbb{R}^n \to \mathbb{R}$ with $P(0) = 0$ 
	is called positive semidefinite (negative semidefinite) on $D$ if
	$P(x) \geqslant 0$ ($P(x) \leqslant 0$) for $x \in D$.
\end{definition}

\begin{definition}
	A continuous function $Q : \left[t_0, \infty\right)\times\mathbb{R}^n \to \mathbb{R}$ 
	with $Q(t, 0) = 0$ is called positive definite (negative definite) on 
	$\left[t_0, \infty\right) \times D$ if there exists a function $\psi \in \mathcal{K}$,
	 such that $\psi\left(\|x\|\right) \leqslant Q(t, x)$ 
	($\psi\left(\|x\|\right) \leqslant -Q(t, x)$) for $t \in \mathbb{T}$, 
	$t \geqslant t_0$, $x \in D$.
\end{definition}

\begin{definition}
	A continuous function $Q : \left[t_0, \infty\right)\times\mathbb{R}^n \to \mathbb{R}$ 
	with $Q(t, 0) = 0$ is called positive semidefinite (negative semidefinite) on 
	$\left[t_0, \infty\right) \times D$ if $Q(t, x) \geqslant 0$ 
	($Q(t, x) \leqslant 0$) for $t \in \mathbb{T}$,	$t \geqslant t_0$, $x \in D$.
\end{definition}

In what follows by $V^\Delta(t,x)$ we denote the full $\Delta$-derivative for 
function $V(x(t))$ along solution of \eqref{general-equation}.

\begin{theorem}
	\label{lyapunov-theorem-1}
	If there exists a continuously differentiable positive-definite function $V$ in a 
	neighborhood of zero with $V^\Delta(t, x)$ negative semidefinite, then the
	equilibrium solution $x = 0$ of equation \eqref{general-equation} is stable.
\end{theorem}

\begin{theorem}
	\label{lyapunov-theorem-2}
	If there exists a continuously differentiable, positive definite function $V$ 
	in a neighborhood of zero and there exists a 
	$\xi \in C_{rd}\left(\left[t_0, \infty\right); \left[0, \infty\right)\right)$ 
	and a $\psi \in \mathcal{K}$, such that
	\[
		V^\Delta(t, x) \leqslant -\xi(t)\psi\left(\left\|x\right\|\right),
	\]
	where
	\begin{equation}
		\label{xi-integral-limit}
		\lim\limits_{t \to \infty} \int\limits_{t_0}^t\xi(s)\,\Delta s = \infty,
	\end{equation}
	then the equilibrium solution $x = 0$ to equation \eqref{general-equation} is asymptotically stable.
\end{theorem}

\begin{theorem}
	\label{lyapunov-theorem-3}
	If there exists a continuously differentiable, positive definite function $V$ 
	in a neighborhood of zero and there exists a 
	$\xi \in C_{rd}\left(\left[t_0, \infty\right); \left[0, \infty\right)\right)$ 
	and a $\psi \in \mathcal{K}$, such that
	\[
		V^\Delta(t, x) \leqslant \xi(t)\psi\left(\left\|x\right\|\right),
	\]
	where \eqref{xi-integral-limit} holds, then the equilibrium solution $x = 0$
	to equation \eqref{general-equation} is unstable.
\end{theorem}

Here and elsewhere we shall use spectral matrix norm as a norm by default:
\[
	\left\|A\right\|_2 = \sqrt{\lambda_{max}\left(A^{^*}A\right)}.
\]

Now we formulate a base model for Hopfield network dynamics and few important
results about stability of its solutions. Indeed, let us consider dynamic
equation of the type
\begin{equation}
	\label{hopfield-equation}
	x^\Delta(t) = -Bx(t) + Ag(x(t)) + J, 
\end{equation}
where $t \in \mathbb{T}$, $\sup{\mathbb{T}} = +\infty$, $x(t) \in \mathbb{R}^n$, 
$A = \left(a_{ij}\right)$, $i,j = \overline{1,n}$, $B = \operatorname{diag}\left(b_i\right)$, 
$b_i > 0, i = \overline{1,n}$, $J = \left(J_1, \dots, J_n\right)^T$, 
$g(x) = \left(g_1(x_1), \dots, g_n(x_n)\right)^T$. Also, $\bar{b} = \max\limits_i\left\{b_i\right\}$,
$\underline{b} = \min\limits_i\left\{b_i\right\}$.
Conceptual meaning of model's components will be clarified below. 

We assume on system \eqref{hopfield-equation} as follows.
\begin{itemize}
	\item[$S_1$.] The vector-function $f(x) = -Bx + Ag(x) + J$ is regressive.
	\item[$S_2$.] There exist positive constants $M_i > 0, i = \overline{1,n}$, 
		such that $\left|g_i(x)\right| \leqslant M_i$ for all $x \in \mathbb{R}$.
	\item[$S_3$.] There exist positive constants $\lambda_i > 0, i = \overline{1,n}$
		such that $\left|g_i(x') - g_i(x'')\right| \leqslant \lambda_i \left|x' - x''\right|$
		for all $x', x'' \in \mathbb{R}$. In what follows we denote 
		$\Lambda = \operatorname{diag}\left(\lambda_i\right)$, $L = \max\limits_i{\lambda_i}$.
\end{itemize}

\begin{definition}
	An $n \times n$ matrix $A$ that can be expressed in the form $A = sE -
	B$, where $E$ is an identity matrix, $B = \left(b_{ij}\right)$ with
	$b_{ij} \geqslant 0, 1 \leqslant i, j \leqslant n$, and $s \geqslant
	\rho(B)$, the maximum of the moduli of the eigenvalues of $B$, is
	called an $M$-matrix.
\end{definition}

It should be noted that $M$-matrix can be characterized in many other ways.
Detailed description of forty such ways one can find in \cite{plemmons1977}.
For our purpose we find useful the following definition.
\begin{definition}
	An $n \times n$ matrix $A$ with non-negative diagonal elements and 
	non-positive off-diagonal ones is called $M$-matrix
	when real part of each eigenvalue of A is positive.
\end{definition}

\begin{lemma}\cite[lemma 5.1.2]{martynyuk2012-green}
	Let assumption $S_3$ be fulfilled. If for every fixed $t \in \mathbb{T}$ 
	the matrix $\left(I - \mu(t)B\right)\Lambda^{-1} - \mu(t)\left|A\right|$
	is an $M$-matrix, the function $f(x) = -Bx + Ag(x) + J$ is regressive.
\end{lemma}

\begin{lemma}\cite[lemma 5.1.1]{martynyuk2012-green}
	If for system \eqref{hopfield-equation} conditions $S_1 - S_3$ are satisfied then there
	exists an equilibrium state $x = x^*$ of system (1) and moreover,
	$\left\|x^*\right\| \leqslant r_0$, where
	\[
		r_0 = \left(
			\sum\limits_{i=1}^n
				\dfrac{1}{b_i^2}
				\left(
					\sum\limits_{j=1}^n{M_j\left|a_{ij}\right|} + \left|J_i\right|
				\right)^2
			\right)^\frac{1}{2}.
	\]
	Besides, if the matrix
	$B\Lambda^{-1} - \left|A\right|$ is an $M$-matrix, this equilibrium
	state is unique.
	\label{lemma-unique-equilibrium}
\end{lemma}

And last result we need is so-called Gershgorin circle theorem.  Let $A$ be a
complex $n\times n$ matrix, with entries $a_{ij}$.  For $i\in \{1,\dots ,n\}$
let $\rho_{i}$ be the sum of the absolute values of the non-diagonal entries in
the $i$-th row.  Let $D(a_{ii}, \rho_{i})$ be the closed disc centered at
$a_{ii}$ with radius $\rho_{i}$. Such a disc is called a Gershgorin disc.

\begin{theorem}[Gershgorin circle theorem]
	Every eigenvalue $\nu$ of $A$ lies within at least one of the Gershgorin discs $D(a_{ii}, \rho_{i})$, i.~e.
	there exists $i \in \left\{1, \dots, n\right\}$ such that
	\[
		\left|\nu - a_{ii}\right| \leqslant \rho_i = \sum\limits_{j \neq i}\left|a_{ij}\right|.
	\]
\end{theorem}

\section{Main results}
\subsection{Node-to-node game setup.}
Let us define a set of nodes $V = \left\{1, 2, \dots, n\right\}$.
Each member of $V$ is interpreted as player in some matrix game with its
neighbors.  This game repeats at time steps, discrete or continuous. Set of
$i$-th node neighbors we denote by $\Omega_i$, $k_i = \left|\Omega_i\right|$.
Here we consider only one type of matrix game is known as Prisoner's
Dilemma. Each node has two strategies: cooperate (C) and defect (D).  Payoff
matrix is illustrated below:
\[
	P = \begin{pmatrix}
		b-c & b \\ -c & 0
	\end{pmatrix}.
\]
Here $b$ is a benefit provided by node to its co-player, $c$ is a cost of
cooperation and hereafter we assume $b > c$. In this case the strategy of
mutual defection is the only Nash equilibrium, while mutual cooperation is more
acceptable social outcome.

Current state of $i$-th node at moment $t$ we denote by $S_i(t) \in \{0,1\}$,
where zero state represent the defection strategy. Easy to show that at the 
instant of time $t$ node $i$ gets total payoff equal to 
$\displaystyle -k_icS_i(t) + \sum\limits_{j \in \Omega_i}{bS_j(t)}$.
This equation remains correct regardless of the nature of time. 
Hence in what follows we assume $t \in \mathbb{T}$, where $\mathbb{T}$ is time scale.

\subsection{Hopfield network setup.}

Assume reaction of each node in network is governed by simple threshold rule:
\[
	S_i(t) = 
		\begin{cases}
			\displaystyle
			0, & \text{if} \quad -k_icS_i(t) + \sum\limits_{j \in \Omega_i}{bS_j(t)} < U_i, \\
			1, & \text{if} \quad -k_icS_i(t) + \sum\limits_{j \in \Omega_i}{bS_j(t)} \geqslant U_i,
		\end{cases}
\]
where $U_i$ is individually payoff threshold for cooperation.
With the aim of using Hopfield neurons model we transform last threshold rule
to the rule with continuous responses. In the end this transformation will lead us to 
dynamical system on time scale modelling asymptotic behavior of network.

Let the state variable $S_i$ for $i$-th ``neuron'' have the range $[0, 1]$ and
be a continuous and strictly increasing function of the total payoff $u_i$.  In
biological terms $S_i$ and $u_i$ are output and input signal of $i$-th
``neuron'' respectively.  Input--output relation we denote by $g_i(u_i)$, so
$S_i(t) = g_i(u_i(t))$ and $u_i(t) = g^{-1}(S_i(t))$. If some node having
non-zero payoff abruptly loses all connections in the network, its behavior may
be described by simple dynamic equation:
\[
	\begin{cases}
		C_i u_i^\Delta(t) = -\dfrac{u_i(t)}{R_i}, \\
		u_i(0) = u_{i0}.
	\end{cases}
\]
By analogy with electrical circuit theory in last equation $C_i$ is called 
capacitance of $i$-th node and $R_i$ is called its resistance. 
Obviously, without communication node's payoff will decay to zero as
individual intention of cooperation do. 

With communication the node gets additional payoff playing with its neighbours, so dynamic equation
becomes as follows:
\[
	C_i u_i^\Delta(t) = -c k_i S_i(t) + \sum\limits_{j \in \Omega_i}{b S_j(t)} - \dfrac{u_i(t)}{R_i}, 
	\qquad i = \overline{1,n},\\
\]
or, using input--output relation and adjacency matrix of network $D = \left(d_{ij}\right), i,j = \overline{1,n}$,
\[
	C_i u_i^\Delta(t) = -c k_i g_i(u_i) + \sum\limits_j{bd_{ij} g_j(u_j)} - \dfrac{u_i}{R_i}, 
	\qquad i = \overline{1,n}.
\]

Let us introduce two matrices:
\[
	A = 	\begin{pmatrix}
		-\dfrac{k_1c}{C_1}   & \dfrac{bd_{12}}{C_1} & \dots & \dfrac{bd_{1n}}{C_1} \\
		\dfrac{bd_{21}}{C_2} & -\dfrac{k_2c}{C_2}   & \dots & \dfrac{bd_{2n}}{C_2} \\
			\dots   & \dots   & \dots & \dots   \\
		\dfrac{bd_{n1}}{C_n} & \dfrac{bd_{n2}}{C_n} & \dots & -\dfrac{k_nc}{C_n}
		\end{pmatrix}, \qquad
	B = 	\begin{pmatrix}
		\displaystyle
			\dfrac{1}{R_1C_1} & 0 & \dots & 0 \\
			0 & \dfrac{1}{R_2C_2} & 0 & \dots\\
			\dots  \\
			0 & 0 & \dots & \dfrac{1}{R_nC_n}
		\end{pmatrix}.
\]
Then we can rewrite equation in vector form:
\begin{equation}
	u^\Delta(t) = -Bu(t) + Ag(u(t)).
\end{equation}

It is easy to see that without significant changes in arguments we can consider
more general model with constant input for every node in network.  By denoting
this input as vector $J = \left(J_i\right)$, we finally get our main equation
as follows:
\begin{equation}
	u^\Delta(t) = -Bu(t) + Ag(u(t)) + J.
	\label{u-hopfield-equation}
\end{equation}

\subsection{Stability condition for network game}
\begin{theorem}
	\label{stability-theorem}
	For system \eqref{u-hopfield-equation} assume that conditions $S_1-S_3$
	are valid and
	\[
		k_i < \frac{\lambda_i}{R_i\left(c + b\right)}, \qquad i = \overline{1,n}.
	\]
	Then there exists unique equilibrium state $u = u^*$ of system \eqref{u-hopfield-equation} and 
	$\left\|u^*\right\| \leqslant r_0$.
\end{theorem}

\begin{proof}
	Let us show $Q = B\Lambda^{-1} - \left|A\right|$ be an $M$-matrix.
	Using inequality for number of neighbors it is easy to derive estimation for 
	real part of eigenvalues $\nu$ of matrix~$Q$. Indeed, by Gershgorin circle theorem
	$\operatorname{Re}{\nu} > 0$ if and only if $Q_{ii} - \rho_i > 0$ for
	all $i=\overline{1,n}$ , where
	\begin{gather*}
		Q_{ii} = \frac{\lambda_i}{R_iC_i} - \frac{k_i c}{C_i},\\
		\rho_i = \sum\limits_{j \neq i}\left|Q_{ij}\right| = \frac{b}{C_i}\sum\limits_{j \neq i} d_{ij} = 
			\frac{k_i b}{C_i}.
	\end{gather*}
	We can write 
	\[
		Q_{ii} - \rho_i = \frac{\lambda_i}{R_iC_i} - \frac{k_i c}{C_i} - \frac{k_i b}{C_i} = 
				\frac{\lambda_i}{R_iC_i} - k_i\frac{c+b}{C_i} > 0
	\]
	and it is obvious that the upper bound for $k_i$ in theorem's statement does guarantee last inequality.
	Hence $Q$ is an $M$-matrix and now to end the proof it remains to apply lemma~\ref{lemma-unique-equilibrium}.		
\end{proof}

\begin{remark}
	It is interesting to notice that existence of unique stable state in network does not depend
	on nodes' ``capacitance''. 
\end{remark}
\begin{remark}
	For every particular node in network ratio $\frac{\lambda_i}{R_i}$ describes
	its potential activity. If $\lambda_i > k_iR_i$ then node's output reaction 
	can conquer with  its overall ``resistance'' to accept neighbor's behaviour.
\end{remark}
\begin{remark}
	Notice that in theorem key role plays overall payoff scale of the game expressed as sum $b + c$.
\end{remark}

It is particularly remarkable that existence of the unique equilibrium state in
network is robust against vanishing of any particular node. Indeed, nodes
fulfil conditions of theorem \ref{stability-theorem} independently. So if there
exists the unique equilibrium state, vanishing any particular node does not
affect on the conditions for all other nodes. On the other hand, a new node may
easily violate conditions of theorem and break the existence of equilibrium.

In two theorems below we use Lyapunov method to formulate sufficient conditions
for asymptotic stability of stable state in network dynamics. Let $u = u^*$ be 
the unique stable state of \eqref{u-hopfield-equation}, i. e. 
\[
	-Bu^* + Ag(u^*) + J = 0.
\]
By introducing new variable $z = u - u^*$ we obtain dynamical system on time scale $\mathbb{T}$
\begin{equation}
	z^\Delta(t) = -Bz(t) + Ah(z(t)),
	\label{z-hopfield-equation}
\end{equation}
where $h(z) = g(z + u^*) - g(u^*)$. If conditions $S_1-S_3$ are valid for system
\eqref{u-hopfield-equation}, it is easy to see that
\begin{itemize}
	\item[$\Sigma_1$.] The vector-function $f(z) = -Bz + Ah(z)$ is regressive.
	\item[$\Sigma_2$.] $\left|h_i(z)\right| \leqslant 2M_i, i = \overline{1,n}$ for all $z \in \mathbb{R}$.
	\item[$\Sigma_3$.] $\left|h_i(z') - h_i(z'')\right| \leqslant \lambda_i \left|z' - z''\right|, i = \overline{1,n}$
		for all $z', z'' \in \mathbb{R}$.
\end{itemize}
Conditions $\Sigma_1$--$\Sigma_3$ guarantee existence and uniqueness for solution of \eqref{z-hopfield-equation}
on $t \in \left[t_0, +\infty\right)$ for any initial values $z(t_0) = z_0$. 

\begin{theorem}[Size-dependent condition]
	\label{asymp-stability-theorem-1}
	Under the conditions of theorem \ref{stability-theorem} assume that 
	$\sup\mathbb{T} = +\infty$ and $\mu(t) \leqslant \mu^*$ for all $t \in \mathbb{T}$. 
	If inequality 
	\begin{equation}
		\label{size-dependent-condition}
		\sqrt{n}\left(b + c\right)\max\limits_{1 \leqslant i \leqslant n}\dfrac{k_i}{C_i} 
		\leqslant 
		\frac{- 1 - \mu^*\bar{b} + \sqrt{1 + 2\mu^*\left(\bar{b} + \underline{b}\right)}}{\mu^*L}
	\end{equation}
	holds, then unique equilibrium state $u = u^*$ of system \eqref{u-hopfield-equation} 
	is uniformly asympto\-tically stable.
\end{theorem}
\begin{proof}
	Clearly, stability of the trivial solution $z = 0$ of \eqref{z-hopfield-equation} is equivalent to
	stability of the stable state $u^*$ of \eqref{u-hopfield-equation}.
	Let us choose the $V(z) = z^Tz$ as a Lyapunov function. It can easily be checked that
	$V(z)$ is positive definite. 
	If $z(t)$ is $\Delta$-differentiable at the moment $t \in \mathbb{T}^\kappa$, 
	the full $\Delta$-derivative of $V(z(t))$ along solution of \eqref{z-hopfield-equation}
	be as follows
	\begin{align*}
		V^\Delta(z(t)) 
			&= \left(z^T(t) z(t)\right)^\Delta = z^T(t)z^\Delta(t) + \left[z^T(t)\right]^\Delta z(\sigma(t)) = \\
			&= z^T(t)z^\Delta(t) + \left[z^T(t)\right]^\Delta \left[z(t) + \mu(t)z^\Delta(t)\right] = \\
			&= 2z^T(t)\left[-Bz(t) + Ah(z(t))\right] + \mu(t)\left\|-Bz(t) + Ah(z(t))\right\|_2^2.
	\end{align*}
	Since $B$ is diagonal matrix with all positive diagonal elements, 
	it follows that maximal eigenvalue of $B$ is $\bar{b} = \max\limits_i\left\{b_i\right\}$ and
	the same one of $-B$ is $\underline{b} = \min\limits_i\left\{b_i\right\}$.
	Using properties of matrix and vector norms and the fact that $\|h(z(t))\| \leqslant L\|z(t)\|$ 
	it's easy to obtain following estimation:
	\begin{align*}
		V^\Delta(z(t)) 
		&\leqslant 
			-2\underline{b}\|z(t)\|^2 + 2\|z(t)\|\left\|A\right\|_2\left\|h(z(t))\right\| + 
			\mu(t)\left(\bar{b}\|z(t)\| + \left\|A\right\|_2\left\|h(z(t))\right\|\right)^2 \\
		&\leqslant		
			-2\underline{b}\|z(t)\|^2 + 2L\left\|A\right\|_2\left\|z(t)\right\|^2 + 
			\mu(t)\left(\bar{b}\|z(t)\| + L\left\|A\right\|_2\left\|z(t)\right\|\right)^2 \\
		&\leqslant		
			-\left(2\underline{b} - 2L\left\|A\right\|_2 - 
			\mu(t)\left(\bar{b} + L\left\|A\right\|_2\right)^2\right) \left\|z(t)\right\|^2.
	\end{align*}
	It is obvious that $\psi\left(\|z\|\right) = \left\|z\right\|^2$ belongs to class $\mathcal{K}$.
	Let us prove that the function 
	\[
		\xi(t) = 2\underline{b} - 2L\left\|A\right\|_2 - \mu(t)\left(\bar{b} + L\left\|A\right\|_2\right)^2
	\]
	under theorems' assumptions belongs to $C_{rd}\left(\left[t_0, \infty\right); \left[0, \infty\right)\right)$ 
	and fulfills condition \eqref{xi-integral-limit}.

	Indeed, we have $\xi(t) \geqslant 2\underline{b} - 2L\left\|A\right\|_2 - \mu^*\left(\bar{b} + L\left\|A\right\|_2\right)^2$.
	Hence, 
	\begin{align*}
		\lim\limits_{t \to \infty} \int\limits_{t_0}^t
			\xi(s)\,\Delta s &= 
		\lim\limits_{t \to \infty} \int\limits_{t_0}^t
		\left(2\underline{b} - 2L\left\|A\right\|_2 - \mu(t)\left(\bar{b} + L\left\|A\right\|_2\right)^2\right)\,\Delta s 
		\geqslant \\
		&\geqslant 
			\left(2\underline{b} - 2L\left\|A\right\|_2 - \mu^*\left(\bar{b} + L\left\|A\right\|_2\right)^2\right)
			\lim\limits_{t \to \infty} \int\limits_{t_0}^t\,\Delta s = \infty. 
	\end{align*}

	Solving quadratic inequality
	\[
		2\underline{b} - 2L\left\|A\right\|_2 - \mu^*\left(\bar{b} + L\left\|A\right\|_2\right)^2 \geqslant 0
	\]
	with respect to $\left\|A\right\|_2$, we get 
	\begin{equation}
		\frac{-1 - \mu^*\bar{b} - \sqrt{1 + 2\mu^*\left(\bar{b} + \underline{b}\right)}}{\mu^*L}	
			\leqslant
		\left\|A\right\|_2
			\leqslant 
		\frac{-1 - \mu^*\bar{b} + \sqrt{1 + 2\mu^*\left(\bar{b} + \underline{b}\right)}}{\mu^*L}.
		\label{A-norm-inequality}
	\end{equation}
	Clearly, by definition of matrix norm left inequality always holds. We have
	\begin{align*}
		\left\|A\right\|_2 \leqslant \sqrt{n}\left\|A\right\|_\infty 
			&= \sqrt{n}\max\limits_{1 \leqslant i \leqslant n}\sum\limits_{j=1}^n\left|a_{ij}\right| = \\
			&= \sqrt{n}\max\limits_{1 \leqslant i \leqslant n}
				\left\{\dfrac{bd_{i1}}{C_i} + \dots + 
				\dfrac{bd_{i,j-1}}{C_i} +
				\dfrac{k_ic}{C_i} + 
				\dfrac{bd_{i,j+1}}{C_i} +
				\dots + \dfrac{bd_{in}}{C_i}\right\} = \\
			&= \sqrt{n}\max\limits_{1 \leqslant i \leqslant n}
				\left\{\dfrac{k_ic}{C_i} + \dfrac{b}{C_i}\sum\limits_{j\neq i}d_{ij}\right\} = \\
			&= \sqrt{n}\max\limits_{1 \leqslant i \leqslant n}
				\left\{\dfrac{k_ic}{C_i} + \dfrac{k_ib}{C_i}\right\} = \\
			&= \sqrt{n}\left(b + c\right)\max\limits_{1 \leqslant i \leqslant n}\dfrac{k_i}{C_i}.
	\end{align*}
	Now it is easy to see that inequality \eqref{size-dependent-condition} 
	guarantees non-negativity of $\xi(t)$ and condition \eqref{xi-integral-limit}.
	To conclude the proof, it remains to use theorem \ref{lyapunov-theorem-2}. 
\end{proof}
\begin{remark}
	We stress that the left side of inequality \eqref{size-dependent-condition}	gathers
	main parameters of network structure (size and maximum relative nodes' degree).
	On the other hand, the right side combines main characteristics of time scale, nodes' 
	inertia and resistance, rate of input-output response.
\end{remark}

Obviously, for any given matrix $A$ inequality \eqref{A-norm-inequality} can be checked directly.
\begin{corollary}
	Under the conditions of theorem \ref{stability-theorem} assume that 
	$\sup\mathbb{T} = +\infty$ and $\mu(t) \leqslant \mu^*$ for all $t \in \mathbb{T}$. 
	If inequality 
	\begin{equation}
		2\underline{b} - 2L\left\|A\right\|_2 - \mu(t)\left(\bar{b} + L\left\|A\right\|_2\right)^2 \geqslant 0
	\end{equation}
	holds, then unique equilibrium state $u = u^*$ of system \eqref{u-hopfield-equation} 
	is uniformly asympto\-tically stable.
\end{corollary}

\begin{theorem}[Size-independent condition]
	\label{asymp-stability-theorem-2}
	Under the conditions of theorem \ref{stability-theorem} assume that 
	$\sup\mathbb{T} = +\infty$ and $\mu(t) \leqslant \mu^*$ for all $t \in \mathbb{T}$. 
	Let $C_*$ denote the minimal ``capacitance'' in the network and 
	$K^*$ denote the largest node's degree:
	\[
		C_* = \min\limits_{1 \leqslant i \leqslant n} C_i \, , 
			\qquad 
		K^* = \max\limits_{1 \leqslant j \leqslant n}{k_j}.
	\]
	If inequality 
	\begin{equation}
		\label{size-independent-condition}
		\left(b + c\right)\frac{K^*}{C_*} 
		\leqslant 
		\frac{- 1 - \mu^*\bar{b} + \sqrt{1 + 2\mu^*\left(\bar{b} + \underline{b}\right)}}{\mu^*L}
	\end{equation}
	holds, then unique equilibrium state $u = u^*$ of system \eqref{u-hopfield-equation} 
	is uniformly asympto\-tically stable.
\end{theorem}
\begin{proof}
	By repeating the same steps as in previous theorem, we obtain
	\[
		\left\|A\right\|_2
			\leqslant 
		\frac{-1 - \mu^*\bar{b} + \sqrt{1 + 2\mu^*\left(\bar{b} + \underline{b}\right)}}{\mu^*L}.
	\]
	Now if we recall matrix norm inequality 
	$\left\|A\right\|_2^2 \leqslant \left\|A\right\|_1 \cdot \left\|A\right\|_\infty$,
	we get 
	\begin{align*}
		\left\|A\right\|_2^2 
			&\leqslant \left\|A\right\|_1 \cdot \left\|A\right\|_\infty = \\
			&= \max\limits_{1 \leqslant j \leqslant n}\sum\limits_{i=1}^n|a_{ij}| \cdot 
				\left(b + c\right)\max\limits_{1 \leqslant i \leqslant n}\dfrac{k_i}{C_i} \leqslant \\
			&= \max\limits_{1 \leqslant j \leqslant n}
				\left\{\dfrac{bd_{1j}}{C_1} + \dots + 
				\dfrac{bd_{j-1,j}}{C_{j-1}} +
				\dfrac{k_jc}{C_j} + 
				\dfrac{bd_{j+1,j}}{C_{j+1}} +
				\dots + \dfrac{bd_{nj}}{C_n}\right\} \cdot 
			   \left(b + c\right)\max\limits_{1 \leqslant i \leqslant n}\dfrac{k_i}{C_i} \leqslant \\
			&= \max\limits_{1 \leqslant i \leqslant n}
				\left\{\dfrac{k_jc}{C_j} + \dfrac{b}{C_*}\sum\limits_{i \neq j}d_{ij} \right\} \cdot 
			   \left(b + c\right)\dfrac{K^*}{C_*} \leqslant \\
			&= \max\limits_{1 \leqslant i \leqslant n}
				\left\{\dfrac{k_jc}{C_*} + \dfrac{bk_j}{C_*}\right\} \cdot 
			   \left(b + c\right)\dfrac{K^*}{C_*} \leqslant \\
			&= \left(b + c\right)\dfrac{K^*}{C_*} \cdot \left(b + c\right)\dfrac{K^*}{C_*} = 
				\left(b + c\right)^2\left(\dfrac{K^*}{C_*}\right)^2.
	\end{align*}
	It is obvious that inequality \eqref{size-independent-condition} 
	guarantees non-negativity of $\xi(t)$ and condition \eqref{xi-integral-limit}.
	To conclude the proof, it remains to use theorem \ref{lyapunov-theorem-2}. 
\end{proof}
\begin{remark}
	For large network, i. e. $n \gg 1$, size-dependent condition \eqref{size-dependent-condition}
	is unlikely to be fulfilled. In the same time condition \eqref{size-independent-condition}
	can be valid regardless of network's size.
\end{remark}

\begin{theorem}[Rate of convergence]
	\label{exp-stability}
	Under the conditions of theorem \ref{stability-theorem} assume that 
	$-\underline{b} \in \mathcal{R}^+$ and $\underline{b} - L\left\|A\right\|_2 > 0$, 
	where $\underline{b} = \min\limits_i\left\{b_i\right\}$. Then
	\begin{itemize}
		\item[1)] solution $z = 0$ of the following system is exponentially stable:
			\begin{equation}
				z^\Delta(t) = -Bz(t), \quad z(t_0) = z_0, \quad t \geqslant t_0 \in \mathbb{T};
				\label{z-linear-equation}
			\end{equation}
		\item[2)] unique equilibrium state $z = 0$ of the system \eqref{z-hopfield-equation} 
			is exponentially stable on $t \geqslant t_0 \in \mathbb{T}$ and the following estimation holds:
			\begin{equation}
				\left\|z(t)\right\| \leqslant  
					\|z_0\| \cdot  e_{-\left(\underline{b} - L\left\|A\right\|_2\right)}(t, t_0). 
				\label{z-exp-estimation}
			\end{equation}
	\end{itemize}
\end{theorem}

\begin{proof}
	Since $B = \operatorname{diag}\left(b_i\right)$ it is
	easy to obtain fundamental matrix $\Phi_{-B}(t, t_0) = \operatorname{diag}
	\left(e_{-b_i}(t, t_0)\right)$. 
	\begin{align*}
		\left\|\Phi_{-B}(t, t_0)\right\| &= 
		\sqrt{\lambda_{\max}\left(\Phi_{-B}^T\Phi_{-B}\right)} = 
		\sqrt{\lambda_{\max}\left(\operatorname{diag}
			\left(e^2_{-b_i}(t, t_0)\right)\right)} = \\
		&=
		\max\limits_{1 \leqslant i \leqslant n} 
			\left|e_{-b_i}(t, t_0)\right| = 
		e_{-\underline{b}}(t, t_0).
	\end{align*}
	It proofs 1) (see \cite[Theorem 2.2]{dacunha2005}). 

	The solution of \eqref{z-hopfield-equation} satisfies the variation of
	constants formula \cite{bohner2012dynamic}
	\[
		z(t) = \Phi_{-B}(t, t_0)z_0 + 
			\int\limits_{t_0}^t 
			\Phi_{-B}(t, \sigma(s))Ah(z(s))\,\Delta s.
	\]
	Hence we have
	\begin{align*}
	\left\|z(t)\right\| 
		&\leqslant 
		\left\|\Phi_{-B}(t, t_0) z_0\right\| + 
		\int\limits_{t_0}^t \left\|
				\Phi_{-B}(t, \sigma(s)) \cdot A h(z(s))
			\right\|\,\Delta s \leqslant \\
		&\leqslant
		e_{-\underline{b}}(t, t_0)\|z_0\| + 
		\int\limits_{t_0}^t \left\|e_{-\underline{b}}(t, \sigma(s)) \right\|
				\cdot \left\|A\right\|_2 \left\|h(z(s))\right\|
			\,\Delta s \leqslant \\
		&\leqslant
		e_{-\underline{b}}(t, t_0)\|z_0\| + 
		\int\limits_{t_0}^t 
			\frac{e_{-\underline{b}}(t, s)}{1 - \underline{b}\mu(s)} 
				\cdot \left\|A\right\|_2 L\|z(s)\|
			\,\Delta s.
	\end{align*}
	Multiplying both sides of inequality by $\dfrac{1}{e_{-\underline{b}}(t, t_0)} > 0$
	(due to $-\underline{b} \in \mathcal{R}^+$) we obtain
	\begin{align*}
	\frac{\left\|z(t)\right\|}{e_{-\underline{b}}(t, t_0)}
		&\leqslant 
		\|z_0\| + 
		\int\limits_{t_0}^t 
			\frac{L\left\|A\right\|_2}{1 - \underline{b}\mu(s)} 
				\cdot  
			\frac{e_{-\underline{b}}(t, s)}{e_{-\underline{b}}(t, t_0)}
				\cdot \|z(s)\|
			\,\Delta s = \\ 
		&= \|z_0\| + 
		\int\limits_{t_0}^t 
			\frac{L\left\|A\right\|_2}{1 - \underline{b}\mu(s)} 
				\cdot  
			\frac{\|z(s)\|}{e_{-\underline{b}}(s, t_0)}\,\Delta s. 
	\end{align*}
	Further, by using Grownall's inequality 
	\[
		\frac{\left\|z(t)\right\|}{e_{-\underline{b}}(t, t_0)} 
		\leqslant 
			\|z_0\| \cdot 
			e_{\frac{L\left\|A\right\|_2}{1 - \underline{b}\mu(s)}}(t, t_0),
	\]
	or
	\[
		\left\|z(t)\right\| \leqslant  \|z_0\| \cdot 
		e_{-\underline{b}\oplus\frac{L\left\|A\right\|_2}{1 - \underline{b}\mu(s)}}(t, t_0) = 
		\|z_0\| \cdot 
		e_{-\left(\underline{b} - L\left\|A\right\|_2\right)}(t, t_0). 
	\]
	By conditions of theorem we have $-\left(\underline{b} -
	L\left\|A\right\|_2\right) \in \mathcal{R}^+$.  Therefore, the last
	estimate means that the solution $z = 0$ of \eqref{z-hopfield-equation}
	is exponentially stable. This completes the proof of theorem.
\end{proof}

Proving Theorems \ref{asymp-stability-theorem-1}, \ref{asymp-stability-theorem-2} we obtain
two variants of majorization for $\|A\|_2$. Both of them can be easily used to find 
direct conditions of exponentially stability expressed in the terms of network's structure and
nodes' internal properties.

\begin{corollary}
	Under the conditions of theorem \ref{stability-theorem} assume that 
	$-\underline{b} \in \mathcal{R}^+$ and 
	\[
		\frac{\underline{b}}{L} > 
			\min\left\{
				\sqrt{n}\left(b + c\right)\max\limits_{1 \leqslant i \leqslant n} \frac{k_i}{C_i}, \,\,
				\left(b + c\right)\frac{K^*}{C_*}
			\right\}.
	\]
	Then unique equilibrium state $z = 0$ of the system \eqref{z-hopfield-equation} 
	is exponentially stable on $t \geqslant t_0 \in \mathbb{T}$ and the estimation 
	\eqref{z-exp-estimation} holds.
\end{corollary}

Obviously, the exponential convergence of the solution for \eqref{z-hopfield-equation} to zero
and the solution $u(t)$ for \eqref{u-hopfield-equation} to the unique equilibrium $u^*$ are the same. 

\section{Discussion.}
In this paper we developed conditions for various types of stability in social
networks governed by {\it Imitation of Success} principle.  There isn't direct,
one-to-one correspondence between considered Hopfield neural network model and
original game-based model.  Hence all obtained results can be considered only
as the base of understanding of opinion propagation in social network.

We limited ourselves to the one type of node-to-node game, Prisoner's Dilemma.
Moreover, arguing we deliberately choose few key elements such as $M$-matrix 
characterization, spectral norm estimation etc. Choosing this elements we were
guided by the aim to obtain simple, fast-checkable and meaningful conditions.

It is an open problem to study network dynamics based on the another interesting 
matrix game types. Perhaps, taking into account network's topology or considering
particular type of time scale one can develop more specific, precise, and useful results.

\bibliographystyle{elsarticle-harv} 
\bibliography{references}

\begin{thebibliography}{7}
\expandafter\ifx\csname natexlab\endcsname\relax\def\natexlab#1{#1}\fi
\expandafter\ifx\csname url\endcsname\relax
  \def\url#1{\texttt{#1}}\fi
\expandafter\ifx\csname urlprefix\endcsname\relax\def\urlprefix{URL }\fi

\bibitem[{Bohner and Peterson(2012)}]{bohner2012dynamic}
Bohner, M., Peterson, A., 2012. Dynamic equations on time scales: An
  introduction with applications. Springer Science \& Business Media.

\bibitem[{DaCunha(2005)}]{dacunha2005}
DaCunha, J.~J., 2005. Stability for time varying linear dynamic systems on time
  scales. J. Comput. Appl. Math. 176~(2), 381--410.

\bibitem[{Hoffacker and Tisdell(2005)}]{hoffacker2005}
Hoffacker, J., Tisdell, C.~C., 2005. Stability and instability for dynamic
  equations on time scales. Computers \& Mathematics with Applications 49~(9),
  1327 -- 1334.

\bibitem[{Hopfield(1984)}]{hopfield1984}
Hopfield, J.~J., 1984. Neurons with graded response have collective
  computational properties like those of two-state neurons. Proc. Natl. Acad.
  Sci. USA 81, 3088--3092.

\bibitem[{Manshadi and Saberi(2011)}]{monshadi2011}
Manshadi, V.~H., Saberi, A., 2011. Prisoner's dilemma on graphs with large
  girth. CoRR abs/1102.1038.
\newline\urlprefix\url{http://arxiv.org/abs/1102.1038}

\bibitem[{Martynyuk(2012)}]{martynyuk2012-green}
Martynyuk, A.~A., 2012. Stability theory of solutions of dynamic equations on
  time scales. Phoenix, Kyiv, in Russian.

\bibitem[{Plemmons(1977)}]{plemmons1977}
Plemmons, R.~J., 1977. M-matrix characterizations. Linear Algebra and its
  Applications 18~(2), 175--188.

\end{thebibliography}

\end{document}